\documentclass[reqno]{amsart}
\usepackage{graphicx}
\usepackage{amssymb}

\begin{document}
\centerline{\LARGE \bf On Lattice Coverings by Simplices}

\vspace{0.5cm}
\centerline{\large Fei Xue and Chuanming Zong}

\address{School of Mathematical Sciences, Peking University, Beijing 100871, P. R. China.}
\email{cmzong@math.pku.edu.cn}

\thanks{{\em 2010 Mathematics Subject Classification}: primary 52C17; secondary 52B10, 52C07}

\vspace{0.5cm}
\centerline{\parbox{11.5cm}{\noindent{\bf Abstract.} By studying the volume of a generalized difference body, this paper presents the first nontrivial lower bound for the lattice covering density by $n$-dimensional simplices.}}

\vspace{1cm}
\section{Introduction}

\bigskip
More than 2,300 years ago, Aristotle (384-322 BCE) claimed that the regular tetrahedra can fill the whole space.
In the modern terms, regular tetrahedra of given size can form a tiling of the three-dimensional Euclidean space $\mathbb{E}^3$. In other words, they can form both a packing and a covering in $\mathbb{E}^3$ simultaneously.
If this were true, both the density of the densest packing by congruent regular tetrahedra and the density of the
thinnest covering of $\mathbb{E}^3$ by congruent regular tetrahedra would be one. Unfortunately, Aristotle is wrong and such a tiling is impossible. Aristotle's mistake was discovered in the fifteenth century by Regiomontanus (see \cite{laga12}). Then, one may ask two natural questions: {\it What is the density of the densest packing by congruent regular tetrahedra} and {\it what is the density of the thinnest covering of $\mathbb{E}^3$ by congruent regular tetrahedra}?

As a part of his 18th mathematical problems, D. Hilbert \cite{hilb01} wrote:\lq\lq {\it I point out the following question, related to the preceding one, and important to number theory and perhaps sometimes useful to physics and chemistry: How can one arrange most densely in space an infinite number of equal solids of given form, e.g., spheres with given radii or regular tetrahedra with given edges $($or in prescribed position$)$, that
is, how can one so fit them together that the ratio of the filled to the unfilled space may be as great as possible}?" Since then, many mathematicians made contributions (mistakes as well) to tetrahedra packings. For the complicated history, we refer to \cite{laga12}.

Covering, in certain sense, is a counterpart of packing. Let $K$ denote a convex body in $\mathbb{E}^n$ and let $C$ denote a centrally symmetric one. In particular, let $B_n$, $T_n$ and $W_n$ denote the $n$-dimensional unit ball, the $n$-dimensional regular simplex with unit edges, and the $n$-dimensional unit cube $\{ {\bf x}:\ 0\le |x_i|\le {1\over 2}\}$, respectively. We call $\mathcal{K}=\{K_i:\ K_i ~\mbox{are congruent to}~ K\}$ a covering of $\mathbb{E}^n$ if
$\bigcup_{K_i\in \mathcal{K}} K_i=\mathbb{E}^n$. For such a $\mathcal{K}$ we define an density
$$\theta (\mathcal{K})= \liminf_{ \ell \to \infty} \frac{ {\rm vol} (\mathcal{K} \cap \ell W_n)}{{\rm vol} (\ell W_n)}.$$
Then, we define the {\em congruent covering density}, the {\em translative covering density} and
the {\em lattice covering density} of $K$ respectively as
$$\theta^c(K) = \min_{\mathcal{K}} \{\theta (\mathcal{K}): \mathcal{K}~ \mbox{ a general covering}\},$$
$$\theta^t(K) = \min_{\mathcal{K}} \{\theta (\mathcal{K}): \mathcal{K} ~\mbox{uses translates of} ~K\}$$
and
$$\theta^l(K) = \min_{\mathcal{K}} \{\theta (\mathcal{K}): \mathcal{K} ~\mbox{is a  lattice covering} \}.$$
In fact, for $\theta^c(K)$, $\theta^t(K)$ and $\theta^l(K)$, the unit cube $W_n$ in the definition of $\theta (\mathcal{K})$ can be replaced by any other fixed convex body. In addition, both $\theta^t(K)$ and $\theta^l(K)$
are invariant under non-singular affine linear transformations. Clearly, for these numbers we have
$$1\le \theta^c(K)\le \theta^t(K)\le \theta^l(K).$$

Let $\Lambda $ be a lattice with determinant ${\rm det}(\Lambda )$, and let $\mathcal{L}$ denote the family of all lattices $\Lambda $ such that $K+\Lambda $ is a covering of $\mathbb{E}^n$. Then $\theta^l(K)$ can be reformulated as
$$\theta^l(K)=\min_{\Lambda\in \mathcal{L}}{{{\rm vol}(K)}\over {{\rm det}(\Lambda )}}.$$

In 1939, Kerschner \cite{kers39} proved
$$\theta^c(B_2)=\theta^t(B_2)= \theta^l(B_2)={{2\pi }\over {\sqrt{27}}}.$$
In 1946 and 1950, L. Fejes T\'oth \cite{feje46} and \cite{feje50} proved that
$$\theta^t(C)=\theta^l(C)\le {{2\pi}\over {\sqrt{27}}}$$
holds for all two-dimensional centrally symmetric convex domains, where equality is attained precisely for the
ellipses. In 1950, F\'ary \cite{fary50} proved that $\theta^l(K)\le 3/2$ holds for all two-dimensional convex domains and the equality holds if and only if $K$ is a triangle. It is trivial that $\theta^c(T_2)=1$. However, the fact $\theta^t(T_2)=3/2$ was proved only in 2010 by Januszewski \cite{janu10}. Even in the plane, the following basic problems are still open (see p.19 of \cite{bras05}):

\medskip\noindent
{\bf Conjecture 1.} {\it For every two-dimensional centrally symmetric convex domain $C$ we have}
$$\theta^c(C)=\theta^l(C).$$

\medskip\noindent
{\bf Conjecture 2.} {\it For every two-dimensional convex domain $K$ we have}
$$\theta^t(K)=\theta^l(K).$$

\medskip
In $\mathbb{E}^3$, our knowledge about $\theta^c(K)$, $\theta^t(K)$ and $\theta^l(K)$ is very limited. In fact, except the five types of parallelohedra $P$ which can tile the whole space and therefore $\theta^c(P)=\theta^t(P)=\theta^l(P)=1$ (see \cite{fedo85}), the only known exact result is
$$\theta^l(B_3)={{5\sqrt{5}\pi }\over {24}}=1.463503\ldots ,$$
which was first established by Bambah \cite{bam54} in 1954 (different proofs were discovered by Barnes \cite{bar56} and Few \cite{few56}). About 2000, a particular lattice tiling was independently discovered by \cite{fidu98} and
\cite{doug04} which implies
$$\theta^l(T_3)\le {{125}\over {63}}.$$
In 2006, Conway and Torquato \cite{conw06} discovered a tetrahedra covering which implies
$$\theta^c(T_3)\le {9\over 8}.$$

In $n$-dimensional space, through the works of Bambah, Coxeter, Davenport, Erd\"os, Few, Watson and in particular Rogers, we know that
$$\theta^t(K)\le n\log n+n\log\log n+5n,$$
$$\theta^l(K)\le n^{\log_2\log_en +c},$$
and
$${n\over {e\sqrt{e}}}\ll \theta^t(B_n)\le \theta^l(B_n)\le c\cdot n(\log_e n)^{{1\over 2}\log_2 2\pi e}.$$

In this paper, we prove the following results:

\medskip\noindent
{\bf Theorem 1.} {\it For any pair of positive numbers $k$ and $m$, we have}
$${{{\rm vol} (kT_n-mT_n)}\over {{\rm vol} (T_n)}}=\sum_{i=0}^n{n\choose i}^2k^im^{n-i}.$$

\medskip\noindent
{\bf Theorem 2.} {\it When $n\ge 3$, we have}
$$\theta^l(T_n)\ge 1+{1\over {2^{3n+7}}}.$$

\vspace{0.5cm}
\section{Generalized Difference Bodies}

\bigskip
In 1904, to study lattice packing of convex bodies, Minkowski \cite{min04} introduced the {\it difference body} $D(K)$ of $K$. Namely,
$$D(K)=\{ {\bf x}_1-{\bf x}_2:\ {\bf x}_i\in K\}.$$
In 1920, Blaschke \cite{blas20} asked for bounds for the volume of $D(K)$ in terms of the volume of $K$. Through the works of Blaschke, Bonnesen, Estermann, Fenchel, Rademacher, S\"uss and in particular the surprising work of Rogers and Shephard \cite{roge57} (also see \cite{rog64}), we have
$$2^n\le {{{\rm vol}(D(K))}\over {{\rm vol} (K)}}\le {2n\choose n},$$
where the lower bound can be attained if and only if $K$ is centrally symmetric, and the upper bound can be attained if and only if $K$ is a simplex.

Let $\lambda$ be a positive number, to generalize Blaschke's problem, it is natural to ask for bounds for
$${{{\rm vol} (K-\lambda K)}\over {{\rm vol} (K)}}.$$
By the {\it Brun-Minkowski inequality} it follows that
$${{{\rm vol} (K-\lambda K)}\over {{\rm vol} (K)}}\ge (1 +\lambda )^n,$$
where the equality holds if and only if $K$ is centrally symmetric. For the upper bounds, it turns out to be challenging.

\medskip\noindent
{\bf Theorem 1.} {\it Let $T_n$ denote an $n$-dimensional simplex, then we have}
$${{{\rm vol} (\mu T_n-\nu T_n)}\over {{\rm vol} (T_n)}}= \sum_{i=0}^n{n\choose i}^2\mu^i\nu^{n-i}.$$

\medskip\noindent
{\bf Proof.} Let $\{{\bf e}_1, {\bf e}_2, \ldots , {\bf e}_n\}$ denote a standard basis of $\mathbb{E}^n$. Let $\sigma $ be a nonsingular linear transformation from $\mathbb{E}^n$ to $\mathbb{E}^n$. For any pair of convex bodies $K_1$ and $K_2$, both contain the origin ${\bf o}$, we have
$$\sigma (K_1+K_2)=\sigma (K_1) +\sigma (K_2).$$ Therefore, without loss of generality, we assume that
$$T_n=\left\{ (x_1, x_2, \ldots , x_n):\ x_i\ge 0,\ \sum x_i\le 1\right\}.$$
In other words, $T_n={\rm conv}\left\{ {\bf o}, {\bf e}_1, {\bf e}_2, \ldots , {\bf e}_n\right\}.$

Let $F_i$ denote an $i$-dimensional face of $T_n$ which contains the origin. Clearly, $$F_i={\rm conv}\left\{ {\bf o}, {\bf e}_{j_1}, {\bf e}_{j_2}, \ldots , {\bf e}_{j_i}\right\}$$ holds for $i$ different base vectors and $T_n$ has ${n\choose i}$ such faces. For convenience, we enumerate all such faces as $F_{i,j}$, where $j=1, 2, \ldots , {n\choose i}$, and denote the $(n-i)$-dimensional face of $T_n$ containing ${\bf o}$ and orthogonal to $F_{i,j}$ by $F^*_{i,j}$.

Then, one can deduce that
$$\mu T_n-\nu T_n=\bigcup_{i=0}^n\ \bigcup_{j=1}^{{n\choose i}}\left(\mu F_{i,j}-\nu F^*_{i,j}\right),$$
$${\rm int}\left(\mu F_{i_1,j_1}-\nu F^*_{i_1,j_1}\right)\bigcap {\rm int}\left(\mu F_{i_2,j_2}-\nu F^*_{i_2,j_2}\right)=\emptyset $$
holds for all $(i_1, j_1)\not=(i_2, j_2),$ and
$${\rm vol} \left(\mu F_{i,j}-\nu F^*_{i,j}\right)={1\over {i!}}\cdot {1\over {(n-i)!}}\cdot \mu^i\nu^{n-i}.$$
Therefore, we have
\begin{align*}
{{{\rm vol} (\mu T_n-\nu T_n)}\over {{\rm vol} (T_n)}}& = n!\sum_{i=0}^n{n\choose i}
{1\over {i!}}\cdot {1\over {(n-i)!}}\cdot \mu^i\nu^{n-i}\\
& =\sum_{i=0}^n{n\choose i}^2\mu^i\nu^{n-i}.
\end{align*}

The theorem is proved.\hfill{$\square$}

\medskip\noindent
{\bf Remark 1.} Rogers and Shephard \cite{roge57} did suggest a mean to compute the volume of $D(T_n)$. Our proof here is different from their argument.

\medskip\noindent
{\bf Conjecture 1.} For every $n$-dimensional convex body $K$ we have
$${{{\rm vol} (\mu K-\nu K)}\over {{\rm vol} (K)}}\le \sum_{i=0}^n{n\choose i}^2\mu ^i\nu^{n-i},$$
where the equality holds if and only if $K$ is a simplex.

\medskip\noindent
{\bf Remark 2.} As a special case of Minkowski's theorem on mixed volumes, for any fixed $n$-dimensional convex body $K$ we have
$${\rm vol}(K-\lambda K)=\sum_{i=0}^n{n\choose i}W_i(K,-K)\cdot \lambda^i,$$
where $W_i(K,-K)$ are constants determined by $K$. It was conjectured by Godbersen \cite{godb38} and Makai jr. \cite{maka74} (see p.412 of Schneider \cite{schn93}) that
$$W_i(K,-K)\le {n\choose i} {\rm vol}(K),$$
where the equality holds if and only if $K$ is a simplex. Clearly, Godbersen and Makai's conjecture implies Conjecture 1.

\vspace{0.5cm}
\section{Lattice Coverings by Simplices}

\bigskip
Assume that $K+\Lambda $ is a lattice covering of $\mathbb{E}^n$. Let $\alpha (K,\Lambda )$ denote its {\it star number} and let $\theta (K,\Lambda )$ denote its density. In other words, $\alpha (K,\Lambda )$ is the number of the lattice points ${\bf u}\in \Lambda \setminus\{{\bf o}\}$ such that $K\cap (K+{\bf u})\not= \emptyset ,$ and $\theta(K,\Lambda )={\rm vol} (K)/{\rm det}(\Lambda ).$

\medskip
To show Theorem 2, we need two basic lemmas. Namely,

\medskip\noindent
{\bf Lemma 1 (Hadwiger \cite{hadw69}, see p.283 of \cite{grub87}).}
{\it Let $K+\Lambda $ be a lattice covering of $\mathbb{E}^n$. Then
we have}
$${{{\rm vol} (2K-K)}\over {{\rm vol}(K)}}\cdot \theta (K,\Lambda )\ge \alpha (K,\Lambda ).$$

\medskip\noindent
{\bf Lemma 2 (Rogers and Shephard \cite{roge57}).} {\it An $n$-dimensional convex body $K$ is a simplex if and only if, for any ${\bf x}\in {\rm int}(D(K))$, the intersection $K\cap (K+{\bf x})$ is positively homothetic to $K$.}

\medskip
Let $K+\Lambda $ be a lattice covering and let $K^j$ denote the subset of $K$ such that every point ${\bf x}\in K^j$ is covered by exact $j$ translates in $K+\Lambda $. We have the following basic result.

\medskip\noindent
{\bf Lemma 3.} {\it If $K+\Lambda $ is a covering of $\mathbb{E}^n$, we have}
$$\theta (K,\Lambda )={{{\rm vol}(K)}\over {\sum {1\over j}
{\rm vol} (K^j)}}={{{\rm vol}(K)}\over {{\rm vol}(K)-\sum {{j-1}\over j}{\rm vol} (K^j)}}.$$

\medskip\noindent
{\bf Proof.} Let $K+\Lambda $ be a lattice covering of $\mathbb{E}^n$ with density $\theta (K,\Lambda )$.
Let $\ell $ be a large positive number, let $\ell W_n$ be a big cube with edge length $\ell $, and let $p(\ell )$ denote the number of the lattice points in $\ell W_n$. Clearly we have
$$\theta (K,\Lambda )=\lim_{\ell\to \infty }{{p(\ell )\cdot {\rm vol}(K)}\over {{\rm vol} (\ell W_n)}}.\eqno(1)$$

Let ${\bf x}$ be a point in $\mathbb{E}^n$ and let ${\bf u}$ be a lattice point. We attach a mass density
$$\delta ({\bf x}, K+{\bf u})={1\over j}$$
to ${\bf x}$ with respect to $K+{\bf u}$ if ${\bf x}\in K+{\bf u}$ and ${\bf x}$ belongs to exact $j$ different translates of $K$ in the lattice covering. If ${\bf x}\not\in K+{\bf u}$, we define $\delta ({\bf x}, K+{\bf u})=0$.
Then the total mass density $\delta ({\bf x})$ at ${\bf x}$ is
$$\delta ({\bf x})=\sum_{{\bf u}\in \Lambda }\delta ({\bf x}, K+{\bf u})=j\cdot {1\over j}=1.$$
Therefore we have
\begin{align*}
{\rm vol}(\ell W_n)& = \int_{\ell W_n}\delta ({\bf x})d{\bf x}\\
& = \int_{\ell W_n}\sum_{{\bf u}\in \Lambda }\delta ({\bf x}, K+{\bf u})d{\bf x}\\
& = \sum_{{\bf u}\in \Lambda}\int_{\ell W_n}\delta ({\bf x}, K+{\bf u})d{\bf x}\\
& = (1+o(1))\cdot p(\ell )\int_{\mathbb{E}^n}\delta ({\bf x}, K)d{\bf x}\\
& = (1+o(1))\cdot p(\ell )\int_K\delta ({\bf x}, K)d{\bf x}\\
& = (1+o(1))\cdot p(\ell )\sum {1\over j}{\rm vol} (K^j). \tag{2}
\end{align*}

By (1) and (2), the lemma follows. \hfill{$\square$}

\medskip\noindent
{\bf Proof of Theorem 2.} For convenience, without loss of generality, we assume that $T_n$ is a regular simplex with unit edges in $\mathbb{E}^n$. We consider two cases.

\medskip
\noindent
{\bf Case 1.} $\alpha (T_n, \Lambda )\ge 2^{3n+1}$.

As a corollary of Theorem 1, we get
$${{{\rm vol}(2T_n-T_n)}\over {{\rm vol}(T_n)}}=\sum_{i=0}^n{n\choose i}^22^i\le 2^n{n\choose {[n/2]}}^2\le 2^{3n}.\eqno (3)$$
Therefore, by Lemma 1 we have
$$\theta (T_n,\Lambda )\ge {{\alpha (T_n, \Lambda )}\over {2^{3n}}}\ge 2.$$

\medskip
\noindent
{\bf Case 2.} $\alpha (T_n, \Lambda )\le 2^{3n+1}$.

Let $\partial (K)$ denote the boundary of $K$, and let $\overline{\rm vol}(X)$ denote the $(n-1)$-dimensional measure of a set $X$ in $\mathbb{E}^n$.

Assume that $T_n$ is intersected by $T_n+{\bf u}_1$, $T_n+{\bf u}_2$, $\ldots $, $T_n+{\bf u}_m$, where $m=\alpha (T_n, \Lambda )$. Then, we have
$$\partial (T_n)=\bigcup_{i=1}^m\left(\partial (T_n)\cap (T_n+{\bf u}_i)\right)$$
and therefore
$$\overline{\rm vol}\left(\partial (T_n)\cap (T_n+{\bf u}_k)\right)\ge {1\over m}\overline{\rm vol}(\partial (T_n))$$
holds at least for one of these translates.

By Lemma 2, we know that $T_n\cap (T_n+{\bf u}_k)$ is homothetic to $T_n$. Assuming that
$$T_n\cap (T_n+{\bf u}_k)=\lambda T_n+{\bf y}$$
holds for some suitable positive number $\lambda $ and a point ${\bf y}$, one can deduce that
$$n\cdot \lambda^{n-1}\cdot {\rm vol}(T_{n-1})\ge {1\over m}\cdot (n+1)\cdot {\rm vol}(T_{n-1}),$$
$$\lambda \ge \left( {{n+1}\over {mn}}\right)^{{1\over {n-1}}},$$
$${\rm vol} \left(T_n\cap (T_n+{\bf u}_k)\right)\ge \left( {{n+1}\over {mn}}\right)^{{n\over {n-1}}}{\rm vol}(T_n),$$
and therefore, when $n\ge 3$,
\begin{align*}
\theta (T_n, \Lambda )& = {{{\rm vol}(T_n)}\over {{\rm vol}(T_n)-\sum {{j-1}\over j}{\rm vol} (T_n^j)}}\\
& \ge {{{\rm vol}(T_n)}\over {{\rm vol}(T_n)-{1\over 2}\sum {\rm vol} (T_n^j)}}\\
& \ge {{{\rm vol}(T_n)}\over {{\rm vol}(T_n)-{1\over 2}{\rm vol} (T_n\cap (T_n+{\bf u}_k))}}\\
& \ge {1\over {1-{1\over 2}\left({{n+1}\over {mn}}\right)^{{n\over {n-1}}}}}\\
& \ge {1\over {1-2^{-{{(3n+1)n}\over {n-1}}-1}}}\\
& \ge 1 + {1\over {2^{3n+7}}}.
\end{align*}

As a conclusion of the two cases, Theorem 2 is proved. \hfill{$\square$}

\medskip\noindent
{\bf Remark 3.} By careful estimation, the lower bound can be further slightly improved.

\vspace{1cm}\noindent
{\bf Acknowledgements.} This work is supported by 973 Programs 2013CB834201 and 2011CB302401, the National Science Foundation of China (No.11071003), and the Chang Jiang Scholars Program of China.

\bibliographystyle{amsplain}

\end{document}